\documentclass[11p, leqno]{article}
\usepackage{amsmath}
\usepackage{amsfonts}
\usepackage{amssymb}
\usepackage{graphicx}
\usepackage{color}
\newtheorem{theorem}{Theorem}[section]
\newtheorem{lemma}{Lemma}[section]
\newtheorem{corollary}{Corollary}[section]
\newtheorem{proposition}{Proposition}[section]

\newtheorem{remark}[theorem]{Remark}

\def\qed{\hfill $\square$}

\newcommand{\RR}{\mathbb{R}}

\def\Ric{\text{Ric}}

\def\R{\Bbb R}

\def\ty{\tilde y}

\def\ty{\tilde{y}}
\def\tz{\tilde{z}}

\def\Ric{\operatorname{Ric}}

\numberwithin{equation}{section}

\begin{document}

\title{\bf \Large Local Aronson-B\'enilan estimates and
\\entropy formulae for porous medium and \\fast diffusion equations
on manifolds}

\author{ Peng Lu, Lei Ni, Juan-Luis V\'azquez and C\'edric Villani}

%\date{6 June 2008}

\maketitle

\begin{abstract}
{In this work we  derive local gradient and Laplacian estimates of the Aronson-B\'enilan
and Li-Yau type for positive solutions of porous medium equations posed on Riemannian manifolds
with a lower Ricci curvature bound. We also prove similar results for some fast diffusion
equations. Inspired by Perelman's work we discover some new entropy formulae for these
equations. }
\end{abstract}

\section{Introduction}

The porous medium equation (PME for short)
\begin{equation}\label{pme}
\partial_t u= \Delta u^m,
\end{equation}
where $m>1$, is a nonlinear version of the classical heat equation
(case $m=1$). For  various values of $m>1$ it has arisen in
different applications to model diffusive phenomena like groundwater
infiltration (Boussinesq's model,  1903, with $m=2$), flow of gas in
porous media (Leibenzon-Muskat model, $m\ge 2$), heat radiation in
plasmas ($m>4$), liquid thin films moving under gravity ($m=4$),
crowd-avoiding population diffusion ($m=2$), and others. The
mathematical theory started in the 1950's and got momentum in recent
decades  as a nonlinear diffusion problem with interesting
geometrical aspects (free boundaries) and peculiar functional
analysis (like  generating a  contraction semigroup in $L^1$ and in
Wasserstein metrics). We refer to the monograph \cite{VaBookPME} for
an account of the rather complete theory concerning existence,
uniqueness, regularity and asymptotic behavior of PME, mostly in the
setting of the Euclidean space and on open subsets of it, as well as the
different applications.

The mathematical treatment of PME can be done in a more or less
unified way for all parameters $m>1$.  Our main estimates below are
only valid for nonnegative solutions, hence we will keep the
restriction $u\ge 0$. This is reasonable from physical grounds since
$u$ represents a density, a concentration, a temperature or a
height in the usual applications. However, re-writing (\ref{pme}) in
the more general form $\partial_t u= \Delta(|u|^{m-1}u)$, solutions
with changing sign can also be considered, but the theory is less
advanced.  It has been proved that for given initial data $u_0\in
L^1(\RR^n)$ with $u_0\ge 0$, there exists a unique continuous weak
solution $u(x,t)\ge 0$ of the initial value problem of (\ref{pme}),
with a number of properties.

Some of the existence, uniqueness and regularity properties hold true for the so-called
fast diffusion equation (FDE), which is equation \eqref{pme} with $m \in (0,1)$.
FDE appears in plasma physics and in geometric flows
such as the Ricci flow on surfaces and
the Yamabe flow. However, there are marked differences between PME and FDE
that justify a separate treatment of FDE, cf. \cite{DK}, \cite{VazSmooth}.
In particular, the qualitative properties of FDE become
increasingly complex for small $m$, far away from $m=1$ (very fast diffusion).

As is typical of nonlinear problems, the mathematical theory of PME
and FDE is based on a priori estimates. In 1979, Aronson and
B\'enilan obtained a celebrated second-order differential inequality
of the form (\cite{AB79})
\begin{equation}
\sum_i \frac{\partial }{\partial x^i} \left (mu^{m-2} \frac{\partial
u}{\partial x^i} \right ) \ge - \frac{\kappa}{t}, \qquad \kappa :=
\frac{n}{n(m-1)+2} \label{eq1.2}
 \end{equation}
which applies to all positive smooth solutions of \eqref{pme}
defined on the whole Euclidean space\footnote{In dimension $n=1$ the
restriction is $m>0$ for general solutions, but we may keep $m>-1$
for so-called
 maximal solutions (\cite{ERV}).} on the condition that
$m > m_c := 1-2/n$. Note that $\sum \frac{\partial }{\partial x^i}
\left (mu^{m-2} \frac{\partial u}{\partial x^i} \right )= \Delta
\left ( \frac{m}{m-1}u^{m-1} \right )$ when $m \neq 1$. Precisely
for the heat equation   formula (\ref{eq1.2}) takes the form
\begin{equation}
\Delta\log u+\dfrac{n}{2t}\geq0.~\rule[-8pt]{0pt}{22pt}
\label{heateqnRnHarnQuant}
\end{equation}
obtained by setting $m=1$. Estimate (\ref{eq1.2}) has turned out to
be a key estimate in the development of the theory of the PME and
FDE posed in the whole Euclidean space. In \cite{AB79} the estimate
was used to prove the existence of initial value problem for PME and
FDE. However, it has turned out difficult to find variants of
(\ref{eq1.2}) that hold for flows posed on open domains,  unless in
dimension $n=1$.

In 1986 Li and Yau studied a heat type flow on manifolds
(\cite{LY86}). Among other things, they proved the following Li-Yau
differential Harnack inequality. If $\left(  M^{n},g\right)$ is a
complete Riemannian manifold with nonnegative Ricci curvature and $u:
M\times\lbrack 0,\infty)\rightarrow \mathbb{R}$ is a positive
solution to the heat equation $\partial_t u=\Delta u$,
then there is a lower bound for $\Delta_g\log u$ that has the
precise sharp form \eqref{heateqnRnHarnQuant}.  This
extends to the manifold setting the Euclidean case mentioned above.
 Li and Yau also proved a local result, which implies \eqref{heateqnRnHarnQuant}
 when $u$ is a global positive solution and $M$ is a complete Riemannian manifold
 with nonnegative Ricci curvature. More precisely, they proved the following theorem.

 \begin{theorem} Let $(M^n, g)$ be a complete Riemannian manifold satisfying Ricci curvature
 $\Ric(M)\ge -K^2$ for some $K \geq 0$.
 Let $B(\mathcal O, 2R)$ be a ball of radius $2R$ centered at $\mathcal{O}$.
 Assume that $u(x, t)$ is a positive smooth solution to the heat equation
 on $B(\mathcal O, 2R) \times [0, \infty )$. Then for any $\alpha>1$
 the following estimate holds on $B(\mathcal O, R)$
 $$
 \sup_{B(\mathcal O, R)}
 \left(\frac{|\nabla u|^2}{u^2} -\alpha \frac{u_t}{u}\right)\le \frac{C\alpha^2}{R^2}
 \left(\frac{\alpha^2}{\alpha^2-1}+KR\right)+\frac{n\alpha^2 K}{2(\alpha -1)}
 +\frac{n\alpha^2}{2t}.
 $$
Here $C$ is a constant only depending on $n$.
 \end{theorem}

The theory of PME and FDE on manifolds has not been considered until
recently. The extension of the Aronson-B\'enilan estimate
(\ref{eq1.2}) to the PME on a closed Riemannian manifold with
nonnegative Ricci curvature was done in the book \cite[Chapter 10]{VaBookPME}.

In this paper we prove an extension of the Aronson-B\'enilan
estimate to the PME flow for all $m>1$ (Theorem \ref{local-li-yau})
and the FDE flow for $m \in (m_c,1)$ (Theorem \ref{fde-local})  on
complete Riemannian manifolds with Ricci curvature bounded below.
 The estimates are of local type,  hence even on
Euclidean space, they give more information. The estimates look much
better when the Ricci curvature is nonnegative.

 Recall that for the positive solution $u:=\frac{e^{-f}}{(4\pi t)^{\frac{n}{2}}}$, such that $u^{\frac{1}{2}}\in W^{1,2}(M)$,
 to the heat equation, it was shown
in \cite{N1} that
\begin{equation}\label{entr-heat}
\frac{d \mathcal{W}}{dt}=-\int_M 2t\left(\left|\nabla_i \nabla_j f-\frac{1}{2t}g_{ij}
\right|^2 +R_{ij}f_if_j\right)u\, d\mu
\end{equation}
 where
 $$
 \mathcal{W}(t)=\int_M \left(t|\nabla f|^2+f-n\right)u\, d\mu.
 $$
Using a basic identity involved in
the proof of Theorem \ref{local-li-yau}, we also obtain entropy formulae in the style of
Perelman (\cite{P}). This new entropy formula is the PME/FDE analogue of \eqref{entr-heat}.

 \noindent  {\sc Organization.} In \S 2 we introduce the main ideas of the
 regularity question for the PME in the Euclidean setting.
 We then pose the problems that have to be addressed.
 \S \ref{sect.local} contains the new local estimate
 for positive solutions of the PME on Riemannian manifolds with Ricci curvature
 bounded below. The estimate admits a version valid for the FDE if $m \in (m_c,1)$,
 which is developed in \S \ref{sec.FDE}.
 Consequences in the
 form of Harnack inequalities for PME and FDE are derived in Sections \ref{sect.local} and
 \ref{sec.FDE}  respectively. \S \ref{sec.ent} introduces and studies the entropies.

\section{Regularity of solutions of PME and Aronson-B\'enilan estimate}

A key idea in the PME theory comes from the observation that we can
write the equation as a diffusion equation for a substance with
density $u(x,t)\ge 0$:
\begin{equation}
 \partial_t u= \nabla\cdot (c(u)\nabla u),
 \end{equation}
We find a case of density-dependent diffusivity, i.\,e.,  $c(u)=mu^{m-1}$,  so that $c$
vanishes at $u=0$; this makes the equation degenerate parabolic. It also implies the
property of finite propagation, appearance of free boundaries, and limited regularity.
Typical solutions with free boundaries are only H\"older continuous in space and time.

The second key idea in the study of the PME is to write the equation
as a  law of mass conservation
\begin{equation}
 \partial_t u = - \nabla\cdot (u  {\bf V}),
 \end{equation}
 which identifies the speed as ${\bf V}=-mu^{m-2}\nabla u$, and this in turn allows to write
${\bf V} $ as a potential flow, ${\bf V} =-\nabla p$. This gives for the potential the
expression $p=mu^{m-1}/(m-1)$. In the application to gases in porous media the potential
is just the pressure and the linear speed-pressure relation is known as Darcy's law.
Historically, the letter $v$ has been used for the pressure instead of $p$,
and we will keep that tradition.
This variable $v$ is crucial in the study of free boundaries and regularity.
Note that the {\sl pressure} $v := mu^{m-1}/(m-1)$ satisfies
\begin{equation}\label{eq.press}
 \partial_t v=(m-1)v\Delta v + |\nabla v|^2.
\end{equation}
Thus, near the level $u=0$ we have the formal approximation $\partial_t v \sim
|\nabla v|^2$, that can be easily identified as movement of the front with speed
$-\nabla v$. It also means that the equation is approximately first-order so that
 we expect the Lipschitz continuity of $v$ near the free boundary.

Now we turn to regularity estimates for PME.
The question of Lipschitz regularity of the
 pressure was solved in one space dimension, $n=1$, by Aronson
  who proved a local estimate for $v_x$ using the Bernstein technique (\cite{Ar69}):
  a bounded solution defined in a cylinder in space-time $Q=[a,b]\times [0,T]$
  has a uniform bound for $|v_x|$ inside the domain, i.\,e., in $Q'=[a',b']
  \times [T',T]$, with $a<a'<b'<b$, $0<T'<T$. B\'enilan proved that
  in that situation $v_t$ is locally bounded in a similar way (\cite{BeKN}).

The extension of such results to dimension $n>1$ fails, even for
globally defined solutions. Indeed, it was shown in 1993 that the
so-called focusing solutions are not Lipschitz continuous at the
focusing point (\cite{AG93}), though wide classes of solutions can
be Lipschitz continuous under special conditions (\cite{CVW87}). And
other kind of pointwise gradient estimates also failed. The problem
of minimal regularity was reduced to proving H\"older regularity,
and this was done around 1980 by Caffarelli and Friedman
(\cite{CafFr79}, \cite{CafFr80}). The proof of these results and the
whole theory of the porous medium equation in several space
dimensions was greatly affected by the existence of special
one-sided estimates that we discuss next.

\medskip

The Aronson-B\'enilan estimate (\ref{eq1.2}) can be
written as
\begin{equation}\label{AB.est}
\Delta v\ge -\frac{\kappa}t, \quad \text{ when }m>m_c,
\end{equation}
 where we define $v=\log u$ for $m=1$. Note that with
this definition $v\le 0 $ for $0<m<1$ so care must be taken in manipulating
inequalities when dealing with fast diffusion.  Using the
pressure equation \eqref{eq.press},  it immediately implies that for
PME with $m>1$
\begin{equation}\label{ABvt.est}
v_t\ge |\nabla v|^2 -\frac{(m-1) \kappa}{t}v.
\end{equation}
so in particular
\begin{equation}\label{ABut}
v_t\ge -(m-1)\kappa v/t, \quad \mbox{ and } \quad u_t\ge -\kappa u/t.
\end{equation}
Other forms of parabolic Harnack inequalities follow from such estimates,
and lead to H\"older regularity statements easily.
These estimates have been used for all kinds of purposes in the theory,
like existence of solutions in optimal classes of data, or asymptotic behavior,
cf. \cite{VaBookPME}.

A  striking property of the Aronson-B\'enilan estimate is the fact
that the constant $\kappa$ is optimal when $m > m_c$. Recall the
Barenblatt solutions, in terms of the pressure, is given by
\begin{equation}
V_C(x,t)=\frac{(C\,t^{2\kappa/n}- \kappa x^2)_+}{2d \,t}, \quad C>0.
\end{equation}
Equality holds in \eqref{AB.est} for $v=V_C$ on the set $\{V_C>0\}$.
When $m=1$ the estimates (\ref{eq1.2}) is optimal since equality
holds in (\ref{heateqnRnHarnQuant}) for the Gaussian kernel. In some
sense the Barenblatt solutions play for the PME a role that the
Gaussian kernel plays for the heat equation.

Many attempts have been made to obtain an extension of the estimate or a suitable
 variant for problems where the PME or the FDE are not posed on the whole Euclidean
 space: this can take the form of boundary value problems in bounded domains of $\RR^d$,
  the PME posed on a Riemannian manifold, or even better, a local estimates valid in
one of the two above settings. A straight extension of the global
estimate to boundary value problems in several dimensions has not
been done. In the case of homogeneous Dirichlet problems a literal
extension is even false, in view of explicit solutions. So there is
a hope for local estimates. In one space dimension, the local
estimate of $v_{xx}$, hence of $v_t$, from below was obtained in
\cite{Vloc}, and the bound has a correction term involving distance
to the boundary. But the method fails for $n>1$ because it uses the
previous knowledge of the local bound of $v_x$. A modified version
of the local estimate will be the first objective of the present
paper.

Another research direction concerns the extension of
Aronson-B\'enilan estimate to other equations, like  the
$p$-Laplacian equation or reaction-diffusion equations. Some work have been done,
 for example, in \cite{EV} for the $p$-Laplacian heat equation on Euclidean space,
  and recently in \cite{KN} for the $p$-Lapacian heat equation and (local and global)
  doubly nonlinear equation on manifolds.

%%%%%%%%%%%%%%%%%%%%%%%%%%%%%%%%%%%%%%%%%%%%%%%%%%%%%%%%%%%
\section{Local Aronson-B\'enilan estimates for the PME}
\label{sect.local}

We proceed now with the new estimates.  Let $u\ge 0$ be a solution
to the Porous Medium Equation \eqref{pme}, $m> 1$, posed on an
$n$-dimensional complete Riemannian manifold $(M^n,g)$. We will
assume at least a local bound from below for the Ricci tensor. The
initial and boundary-value problems for this equation are usually
formulated in terms of weak solutions, or better
 continuous weak solutions \cite{VaBookPME}. Our local estimate
 is closer related to the result of Li-Yau mentioned in the introduction,
 than that of say \cite{SZ}.

We will work with the pressure $v$, which satisfies equation
(\ref{eq.press}).  We see that $\nabla v= m u^{m-2} \nabla u $, and
in the case $m<2$ this equation only makes sense over $u>0$. In
order to avoid this and other regularity difficulties in our
computations we will assume that the solutions are positive and
smooth everywhere. The smoothness property comes from local
boundedness and positivity of $u$ in view of standard non-degenerate
parabolic theory. Application of our  results for general weak
solutions proceeds in a standard way by approximating, and using the
maximum principle and the local compactness of the classes of
solutions involved. We refrain from more details on this issue, cf.
\cite{VaBookPME}.

\medskip

\noindent {\bf \ref{sect.local}.1} Assuming that $u>0$ we introduce
the quantities $y=|\nabla v|^2/{v}$, and $z=v_t/{v}$ and the
differential operator
$$
\mathcal{L}:= \frac{\partial}{\partial t} -(m-1) v\Delta.
$$
We also introduce the differential expression \ $ F_\alpha := \alpha
z-y$. Using equation (\ref{eq.press}) we can write the equivalent
formulae
\begin{equation}
F_\alpha = (m-1)\Delta v+ (\alpha-1)\frac{v_t}{v}=\alpha (m-1)\Delta v+(\alpha-1)
\frac{|\nabla v|^2}{v}. \label{eq F alpha}
\end{equation}
In particular, $F_1=(m-1)\Delta v$. Though our main goal is to
estimate $F_1$, we will use the localization technique
 of Li and Yau  to estimate $F_\alpha $ for $\alpha>1$.
 One reason is that sometimes the estimate of $F_1$
 is not feasible, for instance, when we want to obtain local estimates.

The goal of this subsection is to calculate a formula for $\mathcal{L}(F_\alpha)$.
The following formula is helpful in the calculation:
\begin{equation}\label{help-for1}
\mathcal{L}\left(\frac{f}{g}\right)=\frac{1}{g}\mathcal{L}(f) -\frac{f}{g^2}
\mathcal{L} (g) +
2(m-1)v \left \langle \nabla \left(\frac{f}{g}\right), \nabla \log g \right\rangle.
\end{equation}
Direct calculation shows the following Bochner-type formulae:

\begin{lemma}\label{bochner-lemma} Let $u$ be a positive smooth solution to
(\ref{pme})  on manifold $(M^n,g)$ for some $m>0$, and let
$v:=\frac{m}{m-1}u^{m-1}$ be the pressure\footnote{Recall that when
$m=1$, we interpret $v=\log u$.}. Then we have
\begin{eqnarray}
\mathcal{L} (v_t)&=& 2\langle\nabla v, \nabla v_t\rangle +F_1 v_t,\\
\mathcal{L} (|\nabla v|^2)&=& 2 |\nabla v|^2 F_1 +2\langle \nabla (|
\nabla v|^2), \nabla v\rangle \nonumber \\
&\quad& -2(m-1) vv_{ij}^2 -2(m-1)v R_{ij}v_i v_j.
\end{eqnarray}
Here we use the notations $v_i$ instead of $v_{,i}$ to
denote partial derivatives, and we write $v_{ij}$ to denote the
Hessian tensor $H(v)_{ij}$ of $v$,  while $v_{ij}^2$ denotes the
standard norm of the Hessian. $R_{ij}$ is the Ricci tensor and
$R_{ij}v_iv_j=Ric(\nabla v, \nabla v)$.
\end{lemma}

The following proposition is a generalization of the computation
carried out in Proposition 11.12 of \cite{VaBookPME}.

\begin{proposition} Let $u$ and $v$ be as in Lemma \ref{bochner-lemma}.
Then
\begin{eqnarray}
\mathcal{L}(F_\alpha)&=&2(m-1) v_{ij}^2 +2(m-1)R_{ij} v_i v_j +2m \langle
\nabla F_\alpha, \nabla v\rangle  \nonumber\\
&\quad& +(\alpha -1) \left(\frac{v_t}{v}\right)^2 +F_1^2. \label{bochner-key}
\end{eqnarray}
\end{proposition}
\noindent {\sl Proof.}
Using (\ref{help-for1}) and Lemma \ref{bochner-lemma} we have that
\begin{eqnarray*}
\mathcal{L} \left(\frac{|\nabla v|^2}{v}\right)&=& \frac{1}{v}\left( 2 |\nabla v|^2
F_1 +2\langle \nabla (|
\nabla v|^2), \nabla v\rangle\right) \\
&\quad& -2(m-1) v_{ij}^2 -2(m-1) R_{ij}v_i v_j-\frac{|\nabla v|^4}{v^2}\\
&\quad& +2(m-1)v \left \langle \nabla \left( \frac{|\nabla v|^2}{v} \right ),
\nabla \log v \right \rangle;\\
\mathcal{L}\left(\frac{v_t}{v}\right)&=&\frac{1}{v}\left(2\langle\nabla v,
\nabla v_t\rangle +F_1 v_t\right)\\
&\quad& -\frac{v_t}{v}\frac{|\nabla v|^2}{v}+2(m-1)v \left \langle \nabla \left(
\frac{v_t}{v}\right), \nabla \log v \right \rangle.
\end{eqnarray*}
Putting together gives that
\begin{eqnarray*}
\mathcal{L} (F_\alpha) &=&2(m-1) v \langle \nabla F_\alpha, \nabla \log v\rangle
+ 2(m-1) v_{ij}^2 +2(m-1) R_{ij}v_i v_j\\
&\quad& +\alpha \frac{v_t}{v} F_1 -2\frac{|\nabla v|^2}{v}F_1 -\alpha \frac{v_t}{v}
\frac{|\nabla v|^2}{v}+\frac{|\nabla v|^4}{v^2}\\
&\quad& +\frac{2}{v}\langle \nabla v, \alpha \nabla v_t -\nabla (|\nabla v|^2)\rangle.
\end{eqnarray*}

Using
$$
\langle \nabla v, \nabla (v F_\alpha )\rangle = v\langle \nabla v, \nabla F_\alpha \rangle
+ F_\alpha |\nabla v|^2,
$$
we can rewrite the last term in the above formula for $\mathcal{L} (F_\alpha)$ as
$$
\frac{2}{v}\langle \nabla v, \alpha \nabla v_t -\nabla (|\nabla v|^2)
\rangle=2\langle \nabla v, \nabla F_\alpha\rangle+2F_\alpha \frac{|\nabla v|^2}{v}.
$$
Hence we get
\begin{eqnarray*}
\mathcal{L} (F_\alpha) &=&2m v \langle \nabla F_\alpha, \nabla \log v\rangle
+ 2(m-1) v_{ij}^2 +2(m-1) R_{ij}v_i v_j\\
&\quad& +\alpha \frac{v_t}{v} F_1 -2\frac{|\nabla v|^2}{v}F_1 -\alpha
 \frac{v_t}{v}\frac{|\nabla v|^2}{v}+\frac{|\nabla v|^4}{v^2}+2F_\alpha
 \frac{|\nabla v|^2}{v}.
\end{eqnarray*}
Note the last five terms simplifies as
$$
\alpha z(z-y) -2y(z-y) -\alpha z y +y^2 +2(\alpha z-y)y=(\alpha -1) z^2 +(z-y)^2.
$$
This completes the proof of the proposition.\qed

\medskip

When $\alpha=1$, (\ref{bochner-key}) becomes the following formula in \cite{VaBookPME}:
\begin{equation}\label{bochner-key2}
\mathcal{L} F_1 =2(m-1) v_{ij}^2 +2(m-1)R_{ij} v_i v_j +2m \langle \nabla F_1,
\nabla v\rangle+F_1^2.
\end{equation}
From this, for positive smooth solution $u$ to (\ref{pme}) with $ m > 1$
on a closed Riemannian manifold of dimension $n$
 with nonnegative Ricci curvature, the following  estimate follows easily
 from maximum principle (see Proposition 11.12 of \cite{VaBookPME}):
\begin{equation}\label{ly1}
F_1 \ge- \frac{(m-1)\kappa}{t}, \quad \kappa:= \frac{n}{n(m-1)+2}.
\end{equation}

\medskip

\noindent {\bf \ref{sect.local}.2} Now we prove a new local estimate for
PME on complete manifolds.
We employ the localization technique of Li and Yau (\cite{LY86}, see also
\cite{KN}). Denote by $B({\mathcal O}, R)$ the ball of radius $R>0$ and
centered ${\mathcal O}$ in $(M^n,g)$, and denote by $r(x)$ the distance
function from ${\mathcal O}$ to $x$.

\begin{theorem}\label{local-li-yau} Let $u$ be a positive smooth solution to PME
(\ref{pme}), $m >1$, on the cylinder $Q:=B({\mathcal O}, R)\times [0, T]$.
Let $v$ be the pressure and let $v_{\max}^{R,T}
:= \max_{B({\mathcal O}, R)\times[0,T]} v$.

\noindent  {\rm (1)} Assume that Ricci curvature $\Ric\ge 0$ on $B({\mathcal O}, R)$.
Then, for any $\alpha>1$ we have
\begin{equation}\label{ly-local1}
\frac{|\nabla v|^2}{v}-\alpha\frac{v_t}{v}\le a\alpha^2\left(\frac{1}{t}+
\frac{v_{\max}^{R,T}}{R^2}\left( C_1+ C_2(\alpha)\right)\right)
\end{equation}
 on $Q' :=B({\mathcal O}, R/2)\times[0,T]$. Here, $a: =\frac{n(m-1)
 }{n(m-1)+2}=(m-1)\kappa$, and the positive constants $C_1$ and $C_2(\alpha)$
 depend also on $m$ and $n$.

\noindent  {\rm (2)} Assume that  $\Ric \ge -(n-1) K^2$ on
$B({\mathcal O}, R)$ for some $K \geq 0$.  Then, for any $\alpha>1$,
we have that on $Q'$,
\begin{equation}
\frac{|\nabla v|^2}{v}-\alpha\frac{v_t}{v}\le
a\alpha^2\left(\frac{1}{t}+C_3(\alpha) K^2v_{\max}^{R,T} \right)
  + a\alpha^2\frac{v_{\max}^{R,T}}{R^2}\left(C_2(\alpha)
+ C_1'(KR) \right)\,. \label{ly-local2}
\end{equation} Here, $a$ and
$C_2(\alpha)$ are as before and the  positive constants
$C_3(\alpha)$ and $C_1'(KR)$ depend also on $m$ and $n$.

Acceptable values of the constants are:
\begin{align*}
& C_1: =40(m-1)(n+2), \qquad  && C_2(\alpha) : =\frac{200a\alpha^2
m^2}{ \alpha -1 } \\
& C_3(\alpha) := \frac{(m-1)(n-1)}{\alpha-1},
&& C_1'( KR):=40(m-1)[ 3+(n-1)(1+KR) ].
 \end{align*}
Note that $C'_1( 0)=C_1$.
\end{theorem}

\medskip

\noindent {\sl Proof.} (i) We start with an auxiliary calculation
about suitable cutoff functions. We take a cut-off function
$\eta(x)$ of the form $\eta(x): =\theta \left(r(x)/R\right)$, where
$\theta(t)$ is a smooth monotone function satisfying the following
conditions $\theta (t)\equiv 1$ for $0\le t \le \frac{1}{2}$,
$\theta(t)\equiv 0$ for $t\ge 1$, $(\theta')^2/\theta \le 40$, and
$\theta''\ge -40\theta\ge -40$. 40 is just a convenient number, it
could be optimized. Direct calculation shows  that on $B({\mathcal
O}, R)$
\begin{equation}
\frac{|\nabla \eta|^2}{\eta}\le \displaystyle\frac{40}{R^2},
\end{equation}
and also
\begin{equation}\label{eq Laplacian of phi}
\Delta \eta \ge -\displaystyle\frac{40((n-1)(1+KR)+1)}{R^2}, \mbox{
if }\, \Ric\ge -(n-1)K^2,
\end{equation}
 with the help of the Laplacian comparison theorem. In particular
$\Delta \eta \ge  -40nR^{-2}$ when  $\Ric\ge 0$. We give some
details about deriving (\ref{eq Laplacian of phi}). Note that
\begin{equation*}
\Delta \eta  =\frac{\theta^{\prime\prime}\left\vert \nabla r \right\vert ^{2}
}{R^{2}}+\frac{\theta^{\prime}\Delta r}{R}.
\end{equation*}
By the Laplacian comparison theorem when $\Ric\ge -(n-1) K^2$
\[
\Delta r \leq(n-1) K \coth( K r).
\]
Since $\coth$ is decreasing, and $\theta^{\prime}=0$ when $r(x)<\frac{1}{2}R$, this
implies
\[
\Delta\eta \geq\frac{-40}{R^{2}}-\frac{\sqrt{40}(n-1)}{R} {K}\coth(\frac{1}{2}KR)
\]
where we have used $\left\vert \theta^{\prime}\right\vert \leq \sqrt{40}$.
Using the inequality ${K}\coth(KR)\leq\frac{1}{R}(1+{K}R)$, one has
\begin{equation*}
\Delta \eta \geq\frac{-40}{R^{2}}-\frac{2 \sqrt{40}(n-1)}{R^{2}}(1+{K}R).
\end{equation*}

\medskip

(ii) To obtain the desired estimates,  we apply the operator
$\mathcal{L}$ to the function $ t\eta (-F_\alpha)$, then
apply the maximum principle argument. Note that if $t\eta
(-F_\alpha) \leq 0$ on $Q$, then (\ref{ly-local2}) follows. Now we
assume $\max_{(x,t) \in Q} t\eta (-F_\alpha) >0$. Let $(x_0, t_0)$
be a point where $t\eta (-F_\alpha)$ achieves the positive maximum.
Clearly we have $t_0 >0$ and at $(x_0, t_0)$
\begin{eqnarray*}
\nabla F_\alpha =-\frac{\nabla \eta}{\eta} F_\alpha, \quad
\mathcal{L}(t\eta (-F_\alpha)) \ge 0.
\end{eqnarray*}

All further calculations in this proof will be at $(x_0,t_0)$.
Denote $C_4: =40((n-1)(1+KR)+1)$,  $\ty: =\eta y= \eta \frac{|\nabla
v|^2}{v}$ and $\tz: =\eta z =\eta \frac{v_t}{v}$. Combining
(\ref{bochner-key}) with the above estimates of $\eta$, we have that
when $\Ric\ge -(n-1)K^2$ on $B(\mathcal O,R)$,
\begin{eqnarray*}
0&\le& \eta \mathcal{L}(t\eta (- F_\alpha )) \\
&\le&-t\eta^2\left( 2(m-1) v_{ij}^2 +2(m-1)R_{ij} v_i v_j  \right)
+2mt\eta^2 \langle \nabla (-F_\alpha), \nabla v\rangle-(\alpha -1)t\eta^2 z^2\\
&\, &-t\eta^2 F_1^2 +2t(m-1)v\eta \frac{|\nabla \eta|^2}{\eta} (-F_\alpha)+
\frac{C_4 }{R^2}(m-1)tv\eta (- F_\alpha)+\eta^2 (-F_\alpha) \\
& \le & - \frac{n(m-1)+2}{n(m-1)}\cdot t(\ty- \tz)^2 +2(n-1)(m-1)K^2 t \ty \eta v
+2mt (\ty-\alpha \tz)|\nabla \eta||\nabla v|\\
&\,& -(\alpha -1)t\tz^2+(\ty-\alpha \tz)\left((m-1)\frac{\left(80+C_4
\right)}{R^2}\cdot tv+1\right).
\end{eqnarray*}
In the last inequality above we have used $v_{ij}^2\ge \frac{(\Delta v)^2}{n}$
 and $(m-1)\Delta v=z-y=F_1$.
Now write $C_5: =80+C_4$ and
$$
(\ty-\tz)^2=\frac{1}{\alpha^2}(\ty-\alpha \tz)^2+2\frac{\alpha
-1}{\alpha^2}(\ty-\alpha \tz) \ty +\left(\frac{\alpha-1}{\alpha}\right)^2
\ty^2.
$$
Also note that
$$
2mt (\ty-\alpha \tz)|\nabla \eta|\cdot |\nabla v|\le \frac{40}{R} mt(\ty-\alpha
\tz)\ty^{1/2} v^{1/2}.
$$
Putting these together we have that when $\Ric\ge -(n-1)K^2$ on $B({\mathcal O}, R)$,
\begin{eqnarray}
0&\le& -\frac{t}{a\alpha^2}(\ty-\alpha \tz)^2+t(\ty-\alpha \tz)\left(
-\frac{2(\alpha-1) }{a\alpha^2}\ty +\frac{40m}{R} \cdot \ty^{1/2} v^{1/2}
+(m-1)\frac{C_5}{R^2}v\right) \nonumber \\
&\,& +(\ty-\alpha \tz)-\frac{1}{a}\left(\frac{\alpha-1}{\alpha}\right)^2
t \ty^2+ 2(m-1)(n-1) K^2 t \ty v\eta -(\alpha -1)t\tz^2. \label{eq major pme}
\end{eqnarray}

\medskip

\noindent (1) When $K=0$, using $-Ax^2+Bx \le \frac{B^2}{4A}$, it
follows from (\ref{eq major pme})
\begin{equation}
0\le -\frac{t}{a\alpha^2}(\ty-\alpha \tz)^2+(\ty-\alpha \tz)\left(\frac{tv}{
R^2}\left(\frac{200 a\alpha^2 m^2}{\alpha -1}+(m-1)C_5\right)
+1\right). \label{eq tem use 5}
\end{equation}
This gives the first estimate, (\ref{ly-local1}).

\medskip

\noindent (2) When $K \neq 0$, in (\ref{eq major pme}) we handle the
$(\ty-\alpha \tz)$-term as in (\ref{eq tem use 5}) with $C_6:
=\frac{200 a\alpha^2 m^2}{\alpha -1}+(m-1)C_5$, and use
$$
-\frac{1}{a}\left(\frac{\alpha-1}{\alpha}\right)^2 t \ty^2+
2(m-1)(n-1)tK^2 \ty v\eta \le C_7 tv^2\,,
$$
where
$$
C_7 :=\frac{(m-1)^2(n-1)^2a\alpha^2 K^4}{(\alpha -1)^2}.
$$
Then the above quadratical inequality (\ref{eq major pme}) on $(\ty-\alpha \tz)$
reduces to
$$
0\le -\frac{t}{a\alpha^2}(\ty-\alpha \tz)^2+\left(C_6 \frac{tv}{R^2}
+1\right)(\ty-\alpha \tz) +C_7tv^2.
$$
This implies that
\begin{eqnarray*}
\ty-\alpha \tz &\le& \frac{a\alpha^2}{2} \left(C_6 \frac{v}{R^2}
+\frac{1}{t}+\sqrt{\left(C_6 \frac{v}{R^2}+\frac{1}{t}\right)^2+
4\frac{C_7}{a\alpha^2}v^2 }\right)\\
&\le& a\alpha^2\left(C_6 \frac{v}{R^2}+\frac{1}{t}+\frac{(m-1)(n-1)K^2}{
\alpha-1}v\right).
\end{eqnarray*}
The claimed result follows from this easily.
\qed

\medskip
Using the local estimate one can generalize Proposition 11.12 of
 \cite{VaBookPME} to noncompact complete Riemannian manifolds with
  nonnegative Ricci curvature.

\begin{corollary} \label{cor global pme}
Let $u(x,t)$, $t \in [0,T]$, be a smooth positive solution of the PME (\ref{pme}) with
$m>1$ on a complete manifold $(M^n,g)$.

\noindent  {\rm (1)} If $(M, g)$ has nonnegative Ricci curvature, then (\ref{ly1}) holds
for $t \in (0,T]$,  provided that $v(x, t)=o(r^2(x))$ uniformly in $t \in (0,T]$.

\noindent  {\rm (2)} If the Ricci curvature $\Ric \ge -(n-1) K^2$ on $M$
for some $K \geq 0$ and $v_{\max} := \max_{M \times [0,T]} v < \infty
$, then for any $\alpha >1$
\begin{equation}
\alpha \frac{v_t}{v} -\frac{|\nabla v|^2}{v} \ge
-(m-1) \kappa \alpha^2 \left(\frac{1}{t}+\frac{(m-1)
(n-1)}{\alpha -1} K^2v_{\max} \right). \label{eq for the Hanrack int}
\end{equation}
\end{corollary}
\noindent {\sl Proof.}
(1) Taking $R\to \infty$ and then $\alpha \to 1$ in (\ref{ly-local1}) we have the result.

(2) Taking $R\to \infty$ in (\ref{ly-local2}) we have the result.
\qed

\medskip
Integrating along minimizing geodesic paths of the local estimate, one can obtain the
following Harnack inequality.
Here we just state the most general form. When $K=0$, if we assume $v_{\min} :=
\min_{M\times [0, T]}v >0$, the estimate simplifies by taking $\alpha\to 1$.

\begin{corollary}\label{hk-local} Same notations and assumptions as in
Theorem \ref{local-li-yau}. Denote $v_{\min}^{R/2, T}$ to be
$\min_{B(\mathcal{O}, \frac{R}{2})\times [0, T]} v $.
Assume that  $\Ric \ge -(n-1) K^2$ on $B({\mathcal O}, R)$
for some $K \geq 0$. Then for any $x_1, x_2\in B(\mathcal{O},
\frac{R}{6})$ and $0\le t_1<t_2 \le T$, and any $\alpha >1$
\begin{align*}
& \frac{v(x_2, t_2)}{v(x_1, t_1)} \\
\geq & \left(\frac{t_1}{t_2}\right)^{a \alpha }
\exp\left(-\frac{\alpha d^2(x_1, x_2)}{4v_{\min}^{R/2, T}(t_2-t_1)}-a\alpha(t_2-t_1)
v_{\max}^{R, T}\left( C_3(\alpha) K^2 + \frac{C_2(\alpha)+C_1'(KR)}{R^2}\right)\right).
\end{align*}
where $d(x_1,x_2)$ is the distance and the constants $C_2(\alpha)$, $C_3(\alpha)$
and $C'_1(KP)$ are as in Theorem \ref{local-li-yau}.
\end{corollary}
\noindent {\sl Proof.} For the minimizing geodesic $\gamma(t)$ joining $(x_1, t_1)$
and $(x_2, t_2)$ we have
\begin{eqnarray*}
\log \left(\frac{v(x_2, t_2)}{v(x_1, t_1)}\right)&=& \int_{t_1}^{t_2}\left( \frac{v_t}{v}
+\langle \frac{\nabla v}{v}, \dot{\gamma}\rangle\right)\, ds\\
&\ge &  \int_{t_1}^{t_2} \left(\frac{v_t}{v}-\frac{|\nabla v|^2}{\alpha v}-\frac{\alpha
|\dot{\gamma} |^2}{4v}\right)\, ds.
\end{eqnarray*}
The result follows from the observation that $\gamma(s)$ lies completely inside $B(\mathcal O,
 \frac{R}{2})$ and the estimate in Theorem \ref{local-li-yau}.
\qed

\medskip
A different way manipulating the integration on geodesic path can have the
following consequence of Corollary \ref{cor global pme}.
This estimates are the analogue of the classical one for the positive solutions
to the heat equation (in viewing for the heat equation $v=\log u$).
\begin{corollary}
Same notations and assumptions as in Corollary \ref{cor global pme}.
We further assume that $v_{\max}: = \max_{M \times [0,T]} v < \infty $.
Let $x_{1},x_{2}\in {M}$ and $0 < t_{1}<t_{2} \leq T$.

\noindent  {\rm (1)} If $(M, g)$ has nonnegative Ricci curvature, then
\begin{equation*}
v(x_2, t_2) -v(x_1, t_1) \geq -{(m-1) \kappa v_{\max}}\log\frac{t_2}{t_1}
-\frac{d^2(x_1,x_2)}{4(t_2-t_1)}.
\end{equation*}

\noindent  {\rm (2)} If Ricci curvature $\Ric \ge -(n-1) K^2$
for some $K \geq 0$, then
for any $\alpha >1$
\begin{align*}
& v(x_2, t_2) -v(x_1, t_1) \\
\geq & -{(m-1) \kappa \alpha v_{\max} }\log\frac{t_2}{t_1}
- \frac{ (m-1)^2 (n-1)\kappa \alpha }{\alpha -1} K^2v^2_{\max}(t_2 -t_1 )
-\frac{\alpha d^2(x_1,x_2)}{4 (t_2-t_1)}.
\end{align*}

\end{corollary}

\noindent  {\sl Proof.}
We only prove (2). Let $\gamma(t)$ to
be a constant speed geodesic with $\gamma(t_1)=x_1$ and $\gamma(t_2) =x_2$.
We compute using (\ref{eq for the Hanrack int})
\begin{align*}
& v(x_2, t_2) -v(x_1, t_1)  \geq \int_{t_1}^{t_2} v_t +\left <\nabla v,
\dot{\gamma} \right> dt \\
& \geq  \int_{t_1}^{t_2} \left (
\frac{1}{\alpha} {|\nabla v|^2} -(m-1) \kappa \alpha \left(\frac{1}{t}+\frac{(m-1)
(n-1)}{\alpha -1} K^2v_{\max} \right) v
- \frac{1}{\alpha} {|\nabla v|^2 } -\frac{\alpha}{4}|\dot{\gamma}|^2  \right) dt \\
& \geq -{(m-1) \kappa \alpha v_{\max} }\log\frac{t_2}{t_1}
- \frac{ (m-1)^2 (n-1)\kappa \alpha }{\alpha -1} K^2v^2_{\max}(t_2 -t_1 )
-\frac{ \alpha d^2(x_1,x_2)}{4 (t_2-t_1)}.
\end{align*}
\qed

Theorem \ref{local-li-yau} can be used to give a local lower
estimate for the Laplacian of $v^\beta$ for $\beta>1$.

\begin{corollary} Same assumption and same notations as in Theorem \ref{local-li-yau}.
Let $1<\beta < m/(m-1)$ fixed. Define $\alpha$ by $\frac{
\alpha-1}{\alpha}= (m-1)(\beta-1)$.

\noindent  {\rm (1)} Assume $\Ric \geq 0$ on $B({\mathcal O}, R)$,
then we have on $Q'$
\begin{equation*}
\Delta v^\beta \ge - \kappa \alpha \beta  \left (v_{\max}^{R,T}
\right )^{\beta-1} \left(\frac{1}{t}+
\frac{v_{\max}^{R,T}}{R^2}\left( C_1+ C_2(\alpha)\right)\right)
\end{equation*}

\noindent  {\rm (2)} Assume that  $\Ric \ge -(n-1) K^2$ for some $K
\geq 0$ on $B({\mathcal O}, R)$. Then we have  on $Q'$
\begin{eqnarray*}
& \Delta v^\beta \ge &
- \kappa \alpha \beta  \left (v_{\max}^{R,T} \right )^{\beta-1}
\left(\frac{1}{t}+C_3(\alpha) K^2v_{\max}^{R,T} \right) \\
& & - \kappa \alpha \beta \frac{\left ( v_{\max}^{R,T} \right )^{\beta}}{R^2}
\left(C_2(\alpha) + C_1'(KR) \right) .
\end{eqnarray*}
where the constants $C_2(\alpha)$, $C_3(\alpha)$
and $C'_1(KP)$ are as in Theorem \ref{local-li-yau}.
  \end{corollary}

\noindent  {\sl Proof.}
Clearly $\alpha >1$. We compute
\begin{align*}
\Delta v^{\beta}& =\beta v^{\beta-1}\left ( \Delta v +  (\beta-1) \frac{ |\nabla v|^2}
{v} \right) \\
& = \frac{\beta }{\alpha(m-1)} \cdot v^{\beta-1} \left (\alpha(m-1) \Delta v +
(\alpha-1) \frac{ |\nabla v^2|}{v} \right ) \\
& =\frac{\beta }{\alpha(m-1)} \cdot  v^{\beta-1}F_\alpha.
\end{align*}
The corollary follows easily from Theorem \ref{local-li-yau}.
\qed

%%%%%%%%%%%%%%%%%%%%%%%%%%%%%%%%

\section{Local Aronson-B\'enilan estimates for FDE}
\label{sec.FDE}

The fast diffusion equation, FDE, is equation \eqref{pme} with $m\in
(0,1)$. However, as we have seen, Aronson-B\'enilan estimate for FDE
on Euclidean space holds only in the range $1>m > m_c:=1-
\frac{2}{n}$, where the relevant constant $\kappa = n/(n(m-1)+2)$ is
still a positive number. This is also the range where the Barenblatt
solutions can be written and play similar role as they play in the
theory of PME. Hence the Aronson-B\'enilan estimate holds on the
range as it would be expected. There is another point needs to be
made. Since $m<1$, the pressure $v=\frac{m}{m-1}u^{m-1}$ is negative
and moreover it is an inverse power of $u$. But $(m-1)v$ is still
positive, hence, as shown, the inequalities \eqref{ABvt.est} and
\eqref{ABut} hold for all $m> m_c$.

\medskip

\noindent {\bf \ref{sec.FDE}.1} Let $u$ be a smooth positive
solution to FDE (\ref{pme}) with $ m \in (m_c,1)$ on a closed
Riemannian manifold $(M^n,g)$ with nonnegative Ricci curvature. From
(\ref{bochner-key2}), and since $m-1<0$ we have
\begin{equation*}
\mathcal{L}(F_1)- 2m\langle \nabla F_1, \nabla
v\rangle \leq   F_1^2-2(1-m)v_{ij}^2\leq  - \left(\frac2{(1-m)n}-1\right)F_1^2.
\end{equation*}
The following follows easily from maximum principle
\begin{equation}
F_1 \le - \frac{(m-1)\kappa}{t}, \qquad \mbox{i\,e.,} \quad \Delta
v\ge -\frac{\kappa}t.
\end{equation}
Note the direction in  the $F_1$-inequality differs from (\ref{ly1})
because of $m<1$; on the other side, the $\Delta v$-inequality is
the same as in the case $m>1$.

\medskip

\noindent {\bf \ref{sec.FDE}.2}
Now we prove a new local estimate for FDE with $m\in (m_c,1)$ on complete manifolds.
In this subsection we employ the same notations as in \S 3.2 and use a similar
localization technique as that of Li and Yau. It turns out that this case technically is slightly harder than the previous case. For example we have to make use of the term $(\alpha-1)t\tz^2$ (which was simply dropped before). As we have seen from \S 4.1 for
$m \in (m_c, 1)$ we should estimate $F_\alpha$ from the above (instead of from the below).
For the local estimate, another difference is that for FDE we will estimate $F_{\alpha}$ for  $\alpha <1$
(instead of $\alpha>1$).

Let $(x_0, t_0)$ be a point where function $t\eta F_\alpha $ achieves
the positive maximum. Clearly we have $t_0 >0$ and at $(x_0, t_0)$
\begin{equation*}
\nabla F_\alpha =-\frac{\nabla \eta}{\eta} F_\alpha, \quad\quad
\mathcal{L}(t\eta F_\alpha ) \ge 0.
\end{equation*}

All further calculation in this proof will be at $(x_0,t_0)$.
Combining (\ref{bochner-key}) with the estimates of $\eta$,
we have that when $\Ric\ge -(n-1)K^2$ on $B(\mathcal O,R)$,
\begin{eqnarray*}
0&\le& \eta \mathcal{L}(t\eta F_\alpha) \\
&\le& t\eta^2\left( 2(m-1) v_{ij}^2 +2(m-1)R_{ij} v_i v_j  \right)+2mt\eta^2
\langle \nabla F_\alpha, \nabla v\rangle +(\alpha -1) t\eta^2 z^2\\
&\, & + t\eta^2 F_1^2 +2t(m-1)v\eta \frac{|\nabla \eta|^2}{\eta} F_\alpha
+\frac{C_4}{R^2}(m-1)tv\eta F_\alpha+\eta^2 F_\alpha \\
&\le&  \frac{n(m-1)+2}{n(m-1)} t(\ty- \tz)^2 - 2 (n-1)(m-1)K^2 t \ty \eta v  + 2mt
(\alpha \tz-\ty)|\nabla \eta||\nabla v|\\
&\,& +(\alpha \tz-\ty)\left((m-1)\frac{\left(80+C_4 \right)}{R^2}\cdot tv+1\right)
-(1-\alpha)t\tz^2.
\end{eqnarray*}
Noticing that now we have $(m-1)v>0$ and $a< 0$ for $m\in (m_c, 1)$.
Proceeding as in the proof of (\ref{eq major pme})
we then have that when $\Ric \geq -(n-1)K^2$ on $B(\mathcal O,R)$,

%%% (we also abbreviate it as $\bar{v}_M^{R, T}$,%%%
\begin{align*}
0\le & -\frac{t}{-a\alpha^2}(\alpha \tz-\ty)^2 + t(\alpha \tz-\ty)\left(\frac{2(1-\alpha)
}{-a\alpha^2}(-\ty) +\frac{40m}{R} \cdot (-\ty)^{1/2} (-v)^{1/2}+(m-1)
\frac{C_5}{R^2}v\right) \\
& +(\alpha \tz-\ty) - \frac{t}{-a}\left(\frac{1- \alpha}{\alpha}\right)^2\ty^2- 2 (m-1)
(n-1)K^2 t \ty \eta v -(1-\alpha)t\tz^2.
\end{align*}

Here we can not do what was done in the proof of (\ref{eq major pme})
since $\frac{2(1-\alpha) }{-a\alpha^2}$ the coefficient in front of $-\ty$ ($>0$)
is positive in stead of being negative.
For simplicity let $X=\alpha \tz-\ty>0$, $Y=-\ty >0$, $\beta=1-\alpha>0$, $\gamma=-a>0$.
By writing
$\tz^2=\frac{1}{\alpha^2}(\alpha \tz -\ty+\ty)^2$, and collecting terms we then have
\begin{eqnarray*}
0&\le& -\frac{t}{\gamma \alpha^2}\left((1+\gamma \beta )X^2 -2\beta (1+\gamma)XY
+\beta (\beta+\gamma) Y^2\right)+2(n-1)[(m-1)v] K^2 t Y\\
&\,& +\frac{40mt}{R}X Y^{\frac{1}{2}}(-v)^{\frac{1}{2}}+\frac{tC_5}{R^2}X\left[
(m-1)v\right] +X.
\end{eqnarray*}
Let $\bar{v}^{R,T}_{\max}: =\max_{B({\mathcal O}, R) \times [0,T] } (-v)$.
Using that for any $\epsilon_1, \epsilon_2>0$
$$
\frac{40mt}{R}XY^{\frac{1}{2}}(-v)^{\frac{1}{2}}\le \frac{2\beta \epsilon_1}{\gamma
\alpha^2}tXY+\frac{t\gamma \alpha^2}{2\epsilon_1\beta} \frac{1600m^2}{R^2}\bar{v}_{\max}^{R,
T} X,
$$
and
$$
2(n-1)[(m-1)v] K^2 t Y \le \gamma \alpha^2 \frac{(n-1)^2(1-m)^2 (\bar{v}_{\max}
^{R, T})^2 K^4}{ \epsilon_2^2 \beta^2}t +
\frac{\epsilon^2_2 \beta^2}{\gamma \alpha^2} Y^2t,
$$
we get
\begin{align*}
0\le& -\frac{t}{\gamma \alpha^2}\left((1+\gamma \beta )X^2 -2\beta (1+\gamma+
\epsilon_1)XY +\beta (\beta+\gamma-\beta\epsilon^2_2) Y^2\right)+X \\
& +\frac{t\gamma
\alpha^2}{2\epsilon_1\beta} \frac{1600m^2}{R^2}\bar{v}_{\max}^{R, T} X
+\frac{tC_5}{R^2}X\left[(m-1)v\right]+\gamma \alpha^2 \frac{(n-1)^2(1-m)^2
(v_{\max}^{R, T})^2 K^4}{\beta^2 \epsilon^2_2}t.
\end{align*}

Choosing $\epsilon_2$ small such that
\begin{equation}
\beta+\gamma -\beta \epsilon_2^2>0, \label{epsilon 1 choose}
\end{equation}
then
\begin{equation*}
-2\beta (1+\gamma+
\epsilon_1)XY +\beta (\beta+\gamma-\beta\epsilon^2_2) Y^2
\geq - \beta \frac{(1+\gamma+
\epsilon_1)^2}{\beta +\gamma -\beta\epsilon^2_2} X^2
\end{equation*}
hence we have
\begin{align*}
0 \le &  -\frac{t}{\gamma \alpha^2}\left(1+\gamma \beta -\beta \frac{(1+\gamma+
\epsilon_1)^2}{\beta +\gamma -\beta\epsilon^2_2}\right)X^2 +X
+\frac{t\gamma \alpha^2}{
2\epsilon_1\beta} \frac{1600m^2}{R^2}\bar{v}_{\max}^{R, T} X \\
&+\frac{tC_5}{R^2}X\left ((1-m)\bar{v}_{\max}^{R, T}\right )+\gamma \alpha^2
\frac{(n-1)^2(1-m)^2 (\bar{v}_{\max}^{R, T})^2 K^4}{\beta^2 \epsilon^2_2}t
\end{align*}

Now for any $\beta \in (0,1)$ we can choose small $\epsilon_1>0$ and $\epsilon_2>0$
such that both (\ref{epsilon 1 choose}) and
\begin{equation}
1+\gamma \beta
-\beta \frac{(1+\gamma+\epsilon_1)^2}{\beta +\gamma -\beta\epsilon^2_2}>0
\label{epsilon 2 choose}
\end{equation}
hold. For example, any $(\epsilon_1,\epsilon_2 ) \in
\left (0, \frac{\gamma \alpha^2}{6(1+\gamma)} \right] \times \left (0, \sqrt{\frac{
\gamma \alpha^2}{ 3(1+\gamma)}} \right ]$ works.
Repeating the argument in the proof of Theorem \ref{local-li-yau} after
(\ref{eq major pme}), we have proved the following result for FDE.

\begin{theorem}\label{fde-local}
Let $B({\mathcal O},R)$ be a ball in a complete Riemannian manifold $(M^n,g)$,
and let $u$ be a positive smooth solution to FDE (\ref{pme}) with $m \in (m_c,1)$
on the cylinder $Q=B({\mathcal O}, R)\times [0, T]$.
Let $v$ be the pressure and let $\bar{v}^{R,T}_{\max}=\max_{B({\mathcal O}, R)
\times [0,T] } (-v)$.

\noindent  {\rm (1)} Assume that $\Ric\ge 0$ on $B({\mathcal O},R)$.
Then for any  $0<\alpha<1$, $\epsilon_1>0$ satisfying
\begin{equation}
\bar{C}_1(a, \alpha, \epsilon_1) := 1+(-a)(1-\alpha) -(1-\alpha) \frac{(1
-a+\epsilon_1)^2}{(1-\alpha) -a}>0, \label{epsilon 2A choice}
\end{equation}
we have
\begin{equation}\label{sdf-ly-local1}
\frac{|\nabla v|^2}{v}-\alpha\frac{v_t}{v}\ge - \frac{(-a)\alpha^2}{\bar{C}_1(a, \alpha,
\epsilon_1)}\left(\frac{1}{t}+\frac{\bar{v}^{R,T}_{\max}}{R^2}\left(
 \bar{C}_2(a, \alpha, \epsilon_1)+\bar{C}_3\right)\right)
\end{equation}
on $Q'=B({\mathcal O}, R/2)\times[0,T]$ where
\begin{align*}
& \bar{C}_2(a, \alpha, \epsilon_1) : =1600m^2\frac{(-a) \alpha^2 }{2\epsilon_1
(1-\alpha)} >0 , \\
& \bar{C}_3 : =40(1-m)(n+2)>0.
\end{align*}

\noindent  {\rm (2)} Assume that  $\Ric \ge -(n-1) K^2$ on $B({\mathcal O},R)$ for
some $K \geq 0$. Then for any $0<\alpha<1$ and $\epsilon_1 >0$,
$\epsilon_2>0$ satisfying (\ref{epsilon 1 choose}) and (\ref{epsilon 2 choose})
%% such that $ 1-\alpha -a -(1-\alpha)\epsilon_2^2 >0$ and
%% $$
%% \bar{C}'_1(a, \alpha, \epsilon_1, \epsilon_2)\doteqdot 1+(-a)(1-\alpha) -(1-\alpha)
%% \frac{(1-a+\epsilon_1)^2}{(1-\alpha) -a +(1-\alpha)\epsilon_2^2}>0
%% $$
we have
\begin{align}
\frac{|\nabla v|^2}{v}-\alpha\frac{v_t}{v} \ge &
- \frac{(-a) \alpha^2}{\bar{C}'_1}\left(\frac{1}{t}
+\bar{C}_4(\alpha, \epsilon_2) \sqrt{\bar{C}'_1}K^2 \bar{v}^{R,T}_{\max} \right) \nonumber \\
& - \frac{(-a) \alpha^2}{\bar{C}'_1}\frac{
 \bar{v}^{R,T}_{\max}}{R^2}\left(\bar{C}_2(a, \alpha, \epsilon_1) +\bar{C}_5(
  KR)\right) \label{fde-ly-local2}
\end{align}
on $Q'$. Here
\begin{align*}
& \bar{C}'_1:=\bar{C}'_1(a, \alpha, \epsilon_1, \epsilon_2) := 1+(-a)(1-\alpha) -(1-\alpha)
\frac{(1-a+\epsilon_1)^2}{(1-\alpha) -a -(1-\alpha)\epsilon_2^2} \\
& \bar{C}_4(\alpha, \epsilon_2):= \frac{(n-1)(1-m)}{(1-\alpha) \epsilon_2} \\
& \bar{C}_5( KR):=40(1-m)[ 3+(n-1)(1+KR) ].
\end{align*}
\end{theorem}

\medskip

Noticing that $\bar{C}_1(a, \alpha, \epsilon_1)\to 1$ as $\alpha\to 1$.
Similar to the proof of Corollary \ref{cor global pme},
by taking $R\to \infty$ and then $\alpha\to 1$,
Theorem \ref{fde-local} has the following consequence for global solution of FDE.

\begin{corollary} \label{cor global fde}
Let $u(x,t)$, $t \in [0,T]$, be a smooth positive solution of the FDE (\ref{pme})
with $m \in (m_c,1)$ on a complete manifold $(M^n,g)$.

\noindent  {\rm (1)} If $(M, g)$ has nonnegative Ricci curvature, then

\begin{equation*}
\frac{|\nabla v|^2}{v}-  \frac{v_t}{v}\ge -\frac{(1-m) \kappa}{t}
\end{equation*}
holds for $t \in (0,T]$,  provided that $|v|(x, t)=o(r^2(x))$ uniformly in $t \in (0,T]$.

\noindent  {\rm (2)} If Ricci curvature $\Ric \ge -(n-1) K^2$ on $M$
for some $K \geq 0$ and $\bar{v}_{\max}
:= \max_{M \times [0,T]}(- v) < \infty $,
then for any $\alpha \in (0,1)$ and $\epsilon_2>0$ satisfying
\begin{equation*}
\bar{C}''_1:= \bar{C}''_1(a, \alpha,  \epsilon_2) := 1+(-a)(1-\alpha) -(1-\alpha)
 \frac{(1-a)^2}{(1-\alpha) -a - (1-\alpha)\epsilon_2^2}>0,
\end{equation*}
we have
\begin{equation*}
 \frac{|\nabla v|^2}{v} - \alpha \frac{v_t}{v} \ge
-(1-m) \frac{\kappa  \alpha^2}{\bar{C}''_1} \left(\frac{1}{t}+\frac{(1-m)
(n-1)}{(1-\alpha)\epsilon_2} \sqrt{\bar{C}_1''}K^2 \bar{v}_{\max} \right).
\end{equation*}
\end{corollary}

\medskip
Integrating (\ref{fde-ly-local2}) along minimal geodesic, we obtain
\begin{corollary}\label{hk-local} Same notations and assumptions as in
Theorem \ref{fde-local}. Denote $\bar{v}_{\min}^{R/2, T}$ to be
$\min_{B(\mathcal{O}, \frac{R}{2})\times [0, T]} (-v) $.
Assume that  $\Ric \ge -(n-1) K^2$ on $B({\mathcal O}, R)$
for some $K \geq 0$. Then for any $x_1, x_2\in B(\mathcal{O},
\frac{R}{6})$ and $0\le t_1<t_2 \le T$, and any $\alpha >1$
\begin{align*}
 \frac{- v(x_2, t_2)}{-v(x_1, t_1)} \leq & \left(\frac{t_2}{t_1}\right)^{\frac{(-a)
\alpha}{\bar{C}_1'  } }
 \cdot \exp\left( \frac{\alpha d^2(x_1, x_2)}{4\bar{v}_{\min}^{R/2, T}(t_2-t_1)} \right ) \\
& \cdot \exp \left (\frac{(-a)\alpha}{\bar{C}_1'}
\bar{v}_{\max}^{R, T}\left( \bar{C}_4(\alpha,\epsilon_2) \sqrt{\bar{C}_1'} K^2
+ \frac{\bar{C}_2(a,\alpha,\epsilon_1)+\bar{C}_5(KR)}{R^2}\right)(t_2-t_1)\right).
\end{align*}
where the constants $\bar{C}_1'$, $\bar{C}_2(a,\alpha,\epsilon_1)$, $\bar{C}_4(\alpha,
\epsilon_2)$ and $\bar{C}_5(KP)$ are defined as in Theorem \ref{fde-local}.
\end{corollary}
\noindent {\sl Proof.} For the minimizing geodesic $\gamma(t)$ joining $(x_1, t_1)$
and $(x_2, t_2)$ we have
\begin{eqnarray*}
\log \left(\frac{-v(x_2, t_2)}{-v(x_1, t_1)}\right)&=& \int_{t_1}^{t_2}\left( \frac{v_t}{v}
+\langle \frac{\nabla v}{v}, \dot{\gamma}\rangle\right)\, ds\\
&\le &  \int_{t_1}^{t_2} \left(\frac{v_t}{v} + \frac{|\nabla v|^2}{\alpha (-v)}
 + \frac{\alpha |\dot{\gamma} |^2}{4(-v)}\right)\, ds.
\end{eqnarray*}
The result follows from the observation that $\gamma(s)$ lies completely inside $B(\mathcal O,
 \frac{R}{2})$ and the estimate in Theorem \ref{fde-local}.
\qed

An integral version of Corollary \ref{cor global fde} is
\begin{corollary}
Same notations and assumptions as in Corollary \ref{cor global fde}.
We further assume that $\bar{v}_{\max}  < \infty $.
Let $x_{1},x_{2}\in {M}$ and $0 < t_{1}<t_{2} \leq T$.

\noindent  {\rm (1)} If $(M, g)$ has nonnegative Ricci curvature, then
\begin{equation*}
v(x_2, t_2) -v(x_1, t_1) \geq  -{(1-m) \kappa \bar{v}_{\max}} \log\frac{t_2}{t_1}
-\frac{d^2(x_1,x_2)}{4(t_2-t_1)}.
\end{equation*}

\noindent  {\rm (2)} If Ricci curvature $\Ric \ge -(n-1) K^2$
for some $K \geq 0$,  then
for any $\alpha \in (0,1)$
\begin{align*}
& v(x_2, t_2) -v(x_1, t_1) \\
\geq &  -{(1-m) \frac{\kappa \alpha}{\bar{C}_1''} \bar{v}_{\max} }\log\frac{t_2}{t_1}
- \frac{ (m-1)^2 (n-1)\kappa \alpha }{(1- \alpha)\sqrt{\bar{C}_1''}\epsilon_2 }
K^2\bar{v}^2_{\max}(t_2 -t_1 )
-\frac{\alpha d^2(x_1,x_2)}{4  (t_2-t_1)}.
\end{align*}
\end{corollary}

\medskip

\begin{corollary} Same assumption and same notations as in Theorem \ref{fde-local}.
Let $\bar{\beta} \in (1, \infty)$ fixed. Define $\alpha \in (0,1)$ by $\frac{
\alpha-1}{\alpha}= (m-1)(\bar{\beta}-1)$.

\noindent  {\rm (1)} Assume $\Ric \geq 0$ on $B({\mathcal O}, R)$,
then we have on $Q'$
\begin{equation*}
\Delta (-v)^{\bar{\beta}} \le  \frac{\kappa \alpha \bar{\beta}}{\bar{C}_1(a,\alpha,\epsilon_1)}
\left (\bar{v}_{\max}^{R,T} \right )^{\bar{\beta}-1} \left(\frac{1}{t}+
\frac{\bar{v}_{\max}^{R,T}}{R^2}\left( \bar{C}_2(a,\alpha,\epsilon_1)+ \bar{C}_3
 \right)\right)
\end{equation*}

\noindent  {\rm (2)} Assume that  $\Ric \ge -(n-1) K^2$ on $B({\mathcal O}, R)$
for some $K \geq 0$. Then we have on $Q'$
\begin{align*}
\Delta (-v)^{\bar{\beta} } \le &
 \frac{\kappa \alpha \bar{\beta} }{ \bar{C}_1'}
 \left (\bar{v}_{\max}^{R,T} \right )^{\bar{\beta}-1}
\left(\frac{1}{t}+\bar{C}_4(\alpha,\epsilon_2) \sqrt{\bar{C}_1'} K^2
\bar{v}_{\max}^{R,T} \right) \\
 & + \frac{\kappa \alpha \bar{\beta}}{\bar{C}_1'}
 \frac{\left (\bar{ v}_{\max}^{R,T} \right )^{\bar{\beta}}}{R^2}
\left(\bar{C}_2(a,\alpha,\epsilon_1) + \bar{C}_5(KR) \right) .
\end{align*}
\end{corollary}

%%%%%%%%%%%%%%%%%%%%%%%%%%%%
\section{Entropy formulae}
\label{sec.ent}
In this section we show that the computation in \S \ref{sect.local} and thermodynamic
considerations from Perelman's \cite{P} (see also \cite{N1, N2}) lead to a formula
like Perelman entropy formula for the porous medium equation. The same formula works for the fast diffusion equation and with good concavity for $1-\frac{1}{n}<m<1$.

Assume that $M$ is a compact manifold. First we derive some auxiliary integral formulae.

\begin{lemma} \label{entr-1} Let $u$ be a positive smooth solution of (\ref{pme})
with $m>0$, and let $v$, $F_1$  be as in \S \ref{sect.local}. Then
\begin{eqnarray}
\frac{d}{d t}\int_M vu\, d\mu &=&\int_M F_1 vu\, d\mu= -m\int_M
|\nabla v|^2 u\, d\mu. \label{dt-dis}
\end{eqnarray}
\end{lemma}

\noindent {\sl Proof.} By (\ref{pme}) and (\ref{eq.press}) we have
$\partial_t(vu)=(m-1)vu\Delta v +|\nabla v|^2u+ v\Delta u^m$. Recall
now that $F_1=(m-1)\Delta v$. Then,
$$
\frac{d}{d t}\int_M vu\, d\mu=\int_M F_1 v u\, d\mu +\int_M
\left(|\nabla v|^2 u+v\Delta u^{m}\right)\, d\mu.
$$
Using the identity $\nabla u^m= u\nabla v$, a simple integration by
parts shows that the second term on the right-hand-side is zero. The
first part (\ref{dt-dis}) follows. Integration by parts shows that
\begin{eqnarray*}
\int_M F_1 vu \, d\mu = \int_M (m-1) (\Delta v) vu\, d\mu = \
m\int_M (\Delta v) u^m\, d\mu = -m\int_M |\nabla v|^2 u\, d\mu.
\end{eqnarray*}\qed

\begin{lemma} \label{entr-2}
\begin{eqnarray}
\frac{d}{dt} \int_M F_1 vu\, d\mu &=& 2\int_M
\left((m-1)\left(v_{ij}^2+R_{ij}v_i v_j\right) +F_1^2\right) v u \,
d\mu. \label{dt-energy}
\end{eqnarray}
\end{lemma}

\noindent {\sl Proof.} Using (\ref{pme}), (\ref{eq.press}), and the
formula $\partial_t={\cal L} + (m-1)v\Delta $, we have
\begin{eqnarray*}
\frac{d}{dt} \int_M F_1 vu\, d\mu &=& \int_M \partial_t
F_1\, v u + F_1 \partial_ t (v u)\,d\mu \\
&=&
\int_M ( {\cal L}F_1 )\, v u\,d\mu  + \int_M (m-1)v(\Delta F_1)\,v u\,d\mu\\
&+& \int_M F_1 [(m-1)vu\Delta v +|\nabla v|^2u+ v\Delta u^m]\,d\mu.
\end{eqnarray*}
We will use (\ref{bochner-key2}) to compute the term with $\cal L$,
and we also use $(m-1)\Delta v=F_1$. Then,
\begin{eqnarray*}
\frac{d}{dt} \int_M F_1 vu\, d\mu &=& \int_M \left( (m-1)v(\Delta
F_1) vu +F_1^2 vu +F_1v \Delta u^{m} \right ) d\mu \\
&\quad& +
2(m-1)\int_M \left(v_{ij}^2+R_{ij}v_i v_j\right)v u\, d\mu +\int_M F_1^2 v u\, d\mu\\
&\quad& +2m \int_M \langle \nabla F_1, \nabla v\rangle v u\, d\mu +\int_M F_1|\nabla
 v|^2 u\, d\mu.
\end{eqnarray*}
Using identity $(m-1)\nabla (v^2u)=(2m-1)vu\nabla v$, and
integrating by parts, we get
\begin{eqnarray*}
(m-1)\int_M (\Delta F_1)v^2 u\, d\mu
 = -(2m-1) \int_M \langle \nabla F_1, \nabla v\rangle vu \, d \mu.
 \end{eqnarray*}
Finally, since $\nabla u^m= u\nabla v$ we also have
 \begin{eqnarray*}
\int_M F_1v \Delta u^{m} \, d \mu &=& -\int_M \langle \nabla F_1,
\nabla v\rangle v u \, d \mu -\int_M F_1|\nabla v|^2 u \, d \mu.
\end{eqnarray*}
 Combining these equalities we prove
(\ref{dt-energy}). \qed

\medskip

Recall the constant
$a=\frac{n(m-1)}{n(m-1)+2}=(m-1)\kappa$. We first define
\begin{equation}
\mathcal{N}_{u}(t) := -t^a\int_M vu\, d\mu.
\end{equation} By Lemma
\ref{entr-1}, we have that
\begin{equation}\label{nash-entr}
\frac{d}{ dt} \mathcal{N}_u(t)=-t^a\int_M \left(F_1+\frac{a}{t}\right)vu\,
 d\mu.\end{equation}
Note that the universal estimate (\ref{ly1}) amounts to
$F_1+\frac{a}{t}\ge 0$.
 Now we define the {\bf entropy} to be
\begin{equation}
\mathcal{W}_u(t) := t\frac{d}{dt}\mathcal{N}_u+\mathcal{N}_u.
\end{equation}
Using the last part of \eqref{dt-dis} to compute the term with
$F_1$, we get
$$
\mathcal{W}_u(t)=t^{a+1}\int_M \left(m\frac{|\nabla v|^2}{v}-\frac{a+1}{t}
\right) v u \, d\mu.
$$

\medskip

Now we show the following Perelman type entropy formula for PME. We
put $b=n(m-1)$, so that  $a=\frac{b}{b+2}$.

\begin{theorem}\label{entr-for}
Let $u$ be a positive smooth solution to (\ref{pme}) with $m>0$.
Let $v$ be the pressure and let $\mathcal{W}_u(t)$ be the entropy defined above.
Then
\begin{eqnarray}
\frac{d}{dt} \mathcal{W}_u(t)&=& -2(m-1)t^{a+1}\int_M \left(\left|v_{ij}+\frac{1}{
(b+2)t}g_{ij}\right|^2+R_{ij}v_i v_j\right) v u\, d\mu \nonumber\\
&\quad& -2t^{a+1}\int_M \left(F_1+\frac{a}{t}\right)^2  v u\, d\mu. \label{entropy-pme}
\end{eqnarray}
\end{theorem}

\medskip

\noindent {\sl Proof.}
Note
$$
\frac{d}{dt} \mathcal{W}_u(t) =\frac{d}{dt}\left(t\frac{d}{dt}
\mathcal{N}_u\right)-t^a\int_M \left(F_1+\frac{a}{t}\right)vu\, d\mu.
$$
By Lemma \ref{entr-1} it is easy to see that
\begin{eqnarray*}
\frac{d}{dt}\left(t\frac{d}{dt} \mathcal{N}_u\right)&=&  \frac{d}{dt}
 \left(- t^{a+1}\int_M F_1 v u \, d\mu +a\mathcal{N}_u \right)\\
&=& -2t^{a+1}\left(\int_M \left((m-1)\left(v_{ij}^2+R_{ij}v_i v_j\right)
+F_1^2\right) v u\, d\mu\right)\\
&\quad& -(a+1)t^a \int_M F_1 v u\, d\mu-at^{a}\int_M\left(F_1+\frac{a}{t}
\right) v u\, d\mu.
\end{eqnarray*}
Hence,
\begin{eqnarray*}
\frac{d}{dt} \mathcal{W}_u(t) &=&-2t^{a+1} \int_M \left((m-1)\left(v_{ij}^2
+R_{ij}v_i v_j\right) +F_1^2\right) v u\, d\mu  \\
&\quad& -(a+1)t^{a}\int_M\left(F_1+\frac{a}{t}\right) v u\, d\mu-(a+1)t^a \int_M
 F_1 v u\, d\mu\\
&=&-2t^{a+1}\int_M \left( (m-1)v_{ij}^2 +(m-1)^2(\Delta v)^2+(m-1)\frac{a+1}{t}
\Delta v \right.\\
&\quad& +\left.\frac{a^2+a}{2t^2}\right)v u\, d\mu
-2(m-1)t^{a+1}\int_M  R_{ij}v_i v_j v u\, d\mu.
\end{eqnarray*}
Observing that $a+1=\frac{2(b+1)}{b+2}$, $\frac{a^2+a}{2}
=\frac{b(b+1)}{(b+2)^2}=(m-1)\frac{n(b+1)}{(b+2)^2}$, hence
\begin{align*}
& (m-1)v_{ij}^2 +(m-1)^2(\Delta v)^2+(m-1)\frac{a+1}{t}\Delta v +\frac{a^2+a}{2t^2} \\
 =& (m-1)v_{ij}^2+(m-1)^2 (\Delta v)^2+2(m-1)\frac{b+1}{(b+2)t}\Delta v
+(m-1)\frac{n(b+1)}{(b+2)^2 t^2} \\
= & (m-1)\left(v_{ij}^2+\frac{2}{(b+2)t}\Delta v +\frac{n}{(b+2)^2
t^2}\right)\\
& + \left((m-1)^2(\Delta v)^2 +\frac{2n(m-1)^2}{(b+2)t}\Delta v
+\frac{n^2 (m-1)^2}{(b+2)^2t^2}\right).
\end{align*}
Formula (\ref{entropy-pme}) follows from completing the squares.
\qed

\begin{corollary} \label{cor monot pme entropy}
Let $(M, g)$ be a closed Riemannian manifold with nonnegative
 Ricci curvature. Assume that $u$ is a positive smooth solution to PME (\ref{pme})
 with $m>1$.
 Then

 \noindent  {\rm (1)} $\frac{d}{dt}\mathcal{N}_u(t)\le 0$ and
$\frac{d}{dt}\mathcal{W}_u(t)\le 0$. In particular
  $\mathcal{N}_u(t)$ is a monotone non-decreasing
 concave function in $\frac{1}{t}$.

 \noindent  {\rm (2)} Any ancient positive solution to (\ref{pme})
 must be a constant.
\end{corollary}
\noindent {\sl Proof.} We only need to justify (2).
For the ancient solution, (\ref{ly1}) implies $F_1\ge 0$. On the other hand,
\begin{equation*}
\int_M F_1 uv\, d\mu =-m\int_M |\nabla v|^2 u\, d\mu
\end{equation*}
which implies that $|\nabla v|=0$.
\qed

\begin{remark} The result can be proved for complete noncompact Riemannian manifolds
with the help of the gradient estimates from the previous section.
The interested readers can find the details for the linear heat equation
case in \cite{Chowetc}.
\end{remark}

In \cite{CT}, the following logarithmic Sobolev type inequality related to PME was proved on $\R^n$. (See also \cite{DD} for Sobolev type inequalities related to FDE.)

{\it For any $m>1$, $f\in L^1(\R^n)\cap L^m (\R^n)$, $|\nabla f^{m-
\frac{1}{2}}|\in L^2(\R^n)$, we have
$$
\left(n+\frac{1}{m-1}\right) \int_{\R^n} f^m  dx\le
\frac{1}{2}\left(\frac{2m}{2m-1}\right)^2 \int_{\R^n} \left|\nabla f^{m
-\frac{1}{2}}\right|^2 dx +A_m(\|f\|_1),
$$
where $A_m(K) := \int_{\R^n} \left(\frac{|x|^2}{2}u_{\infty}+\frac{1}
{m-1}u_\infty^m\right)\, dx$ with $u_\infty := (C-\frac{m-1}{2m} |x|
^2)_{+}^{\frac{1}{m-1}}$ being the Barenblatt-Pattle solution  of order
$m$ and mass $K := \int_{\R^n} u_{\infty}\, dx$.
}

The previous defined entropy $\mathcal{W}_u$  is related to the above Sobolev
 inequality in the following way.
First notice that $A(m, n) := A_m(1)$ is a constant depending only on
the dimension $n$ and $m$. Direct calculation shows that Carrillo-Toscani's
logarithmic Sobolev inequality amounts to that
$$
\mathcal{W}_u\left(\frac{1}{b+2}\right)\ge -2m \left(\frac{1}{b+2}\right)^{a+1}A(m, n).
$$

\medskip

For the fast diffusion equation, the monotonicity of entropy $\mathcal{W}_u$ is
similar to that of PME given in Corollary \ref{cor monot pme entropy}.

\begin{corollary} Let $(M, g)$ be a closed Riemannian manifold with nonnegative
Ricci curvature. Assume that $u$ is a positive smooth solution to FDE (\ref{pme})
with $m <1$. Then

 \noindent  {\rm (1)}  $\frac{d}{dt}\mathcal{N}_u(t)\le 0$ for $m\in (m_c, 1)$.

 \noindent  {\rm (2)}
$\frac{d}{dt}\mathcal{W}_u(t)\le 0$ for $m\in [m_c', 1)$
with $m_c'=1-\frac{1}{n}$.
In particular $\mathcal{N}_u(t)$ is a monotone non-decreasing
 concave function in $\frac{1}{t}$ when $m\in [m_c', 1)$.

  \noindent  {\rm (3)} Any ancient positive solution to (\ref{pme}) with $m \in (m_c,1)$
 must be a constant.
\end{corollary}

\noindent {\sl Proof.} (1) In (\ref{nash-entr}), $v <0$ and by Corollary
\ref{cor global fde}(1) $F_1+\frac{a}{t} \leq 0$.

(2) Notice that
\begin{eqnarray*}-(m-1)\left|v_{ij}+\frac{1}{(b+2)t}g_{ij}\right|^2&\ge&
-\frac{1}{n(m-1)}\left((m-1)\Delta v +\frac{n(m-1)}{(b+2)t}\right)^2\\
&=&-\frac{1}{n(m-1)}\left(F_1+\frac{a}{t}\right)^2.
\end{eqnarray*}
Note $v <0$, the result then follows from (\ref{entropy-pme}) and
$1+(m-1)n\ge 0$ for $m\in [m_c', 1)$.
\qed

\

\section*{Acknowledgments}
{P.L. thanks Institut Henri Poincar\'{e}, where  part of
 this joint work was done,  for the hospitality. L.N. was
supported in part by NSF grant DMS-0504792 and an Alfred P. Sloan
Fellowship, USA. J.L.V. is partially funded by Spanish Project
MTM2005-08760-C02-01. He wishes to thank the Dept. of Mathematics of
the University of California at San Diego,  where he
visited to work on this project. C.V. thanks
J.A.~Carrillo for inspiring conversations.}

\bigskip

%%%%%%%%%%%%%%%%%%%%%%%%%%%%%%%%%%%%%%%%%%%%%%%%%%%%%%%%%%%%%%%%%%%%%%%%%%%%%%%
\bibliographystyle{amsalpha}

\

\begin{center} {\bf Keywords:} {Porous medium equation, fast diffusion equation,
Aronson-B\'enilan estimate, Li-Yau type estimate,  local gradient
bound, flow on manifold, entropy formula. \hfill}
\end{center}

\

{\sc Addresses}

\medskip
{\sc Peng Lu},
Department of Mathematics, University of Oregon,
Eugene, OR 97403, USA

e-mail: penglu@uoegon.edu

\medskip

{\sc Lei Ni},
Department of Mathematics, University of California at San Diego,
La Jolla, CA 92093, USA

e-mail: lni@math.ucsd.edu

\medskip

{\sc Juan Luis V\'azquez},
Departamento de Matem\'{a}ticas, Universidad
Aut\'{o}noma de Madrid,
Campus de Cantoblanco, 28049 Madrid, Spain

e-mail: juanluis.vazquez@uam.es

{\sc Cedric Villani},
Unite de Mathematiques Pures et Appliquees
Ecole normale superieure de Lyon
46 allee d'Italie
F-69364 Lyon Cedex 07, France

e-mail: cvillani@umpa.ens-lyon.fr

\end{document}